\documentclass[12pt]{article}

\usepackage{amsmath}
\usepackage{amssymb}

\newcommand{\C}{{\mathbb C}}
\newcommand{\R}{{\mathbb R}}
\newcommand{\T}{{\mathbb T}}
\newcommand{\D}{{\mathbb D}}
\newcommand{\LL}{{\cal L^{++}}}
\newcommand{\cA}{{\cal A}}
\newcommand{\LM}{{\cal L}^-}
\newcommand{\omz}{\omega(0,E,\Omega)}
\renewcommand{\Im}{{\rm Im}\,}

\newtheorem{theorem}{Theorem}
\newtheorem{theo}{Theorem}
\newtheorem{corollary}[theorem]{Corollary}
\newtheorem{lemma}[theorem]{Lemma}
\newtheorem{definition}[theorem]{Definition}
\newtheorem{proposition}{Proposition}

\title{Classical and new $\log\log$-theorems}
\author{Alexander Rashkovskii}
\date{}

\begin{document}

\maketitle

\begin{abstract} We present a unified approach to celebrated $\log\log$-theorems
of Carleman, Wolf, Levinson, Sj{\"o}berg, Matsaev on majorants of
analytic functions. Moreover, we obtain stronger results by replacing original pointwise bounds with integral ones. The main ingredient is a complete description for
radial projections of harmonic measures of strictly star-shaped
domains in the plane, which, in particular, explains where the $\log\log$-conditions come from.
\end{abstract}

\section{Introduction. Statement of results}

Our starting point is classical theorems due to Carleman, Wolf,
Levinson, and Sj{\"o}berg, on majorants of analytic functions.

\begin{definition} A nonnegative measurable function $M$ on a segment
$[a,b]\subset\R$ belongs to the class $\LL[a,b]$ if
$$
\int_a^b \log^+\log^+ M(t)\,dt<\infty.
$$
\end{definition}
(For any real-valued function $h$, we write $h^+=\max\{h,0\}$,
$h^-=h^+-h$.)

Carleman was the first who remarked a special role
of functions of the class $\LL$ in complex analysis, by proving the following variant of the Liouville theorem.

\medskip
\noindent {\bf Theorem A} (T. Carleman \cite{C}) {\em If an entire
function $f$ in the complex plane $\C$ has the bound
\begin{equation}
\label{eq:A}
|f(re^{i\theta})|\le M(\theta)\quad \forall\theta\in [0,2\pi],\
\forall r\ge r_0,
\end{equation}
 with $M\in \LL[0,2\pi]$, then $f\equiv const$.}

\medskip
This phenomenon appears also in the Phragm\'en--Lindel\"of setting.

\medskip\noindent
{\bf Theorem B} (F. Wolf \cite{W1})
{\em If a holomorphic function $f$
in the upper half-plane $\C_+=\{z\in\C:\:{\rm Im}\,z>0\}$ satisfies the condition
$$
\limsup_{z\to x_0} |f(z)|\le 1\quad\forall x_0\in\R
$$
and for any $ \epsilon>0$ and all $r>R(\epsilon)$, $\theta\in
(0,\pi)$, one has $$|f(re^{i\theta})|\le [M(\theta)]^{\epsilon r}$$
with $M\in \LL[0,\pi]$, then $|f(z)|\le 1 $ on $\C_+$.}

\medskip

The most famous statement of this type is the following
local result known as the Levinson--Sj{\"o}berg theorem.

\medskip
\noindent {\bf Theorem C} (N. Levinson \cite{L}, N. Sj{\"o}berg \cite{S}, F. Wolf \cite{W2}) {\em If a
holomorphic function $f$ in the domain $Q=\{x+iy:\: |x|<1,\ |y|<1\}$
has the bound $$|f(x+iy)|\le M(y)\quad\forall x+iy\in Q,$$ with $M\in \LL[-1,1]$, then for any
compact subset $K$ of $Q$ there is a constant $C_K$, independent of
the function $f$, such that $|f(z)|\le C_K$ in $K$.}

\medskip
For further developments of Theorem C, including higher dimensional
variants, see \cite{D1}, \cite{D2}, \cite{Dy2}, \cite{Dy3},
\cite{G}. Theorems A and B were extended to subharmonic functions in
higher dimensions in \cite{Y}.

A similar feature of majorants from the class $\LL$ was discovered
by Beurling in a problem of extension of analytic functions
\cite{Be}. It also appears in relation to holomorphic functions from
the MacLane class in the unit disk \cite{H}, \cite{Mac}, and in a
description of non-quasi-analytic Carleman classes \cite{Dy1}.

The next result, due to Matsaev, does not look like a
$\log\log$-theorem, however (as will be seen from our
considerations) it is also about the class $\LL$; further results in this direction can be found in \cite{MM}.

\medskip
\noindent
{\bf Theorem D} (V.I. Matsaev \cite{M})
{\em If an entire function $f$ satisfies the relation
$$
\log|f(re^{i\theta})|\ge -Cr^\alpha|\sin\theta|^{-k} \quad \forall\theta\in(0,\pi),\ \forall r>0,
$$
with some $C>0$, $\alpha>1$, and $k\ge 0$, then it has at most normal type
with respect to the order $\alpha$, that is, $\log|f(re^{i\theta})|\le
Ar^\alpha+B$.}

\medskip

All these theorems can be formulated in terms of subharmonic
functions (by taking $u(z)=\log|f(z)|$ as a pattern), however our
main goal is to replace the pointwise bounds like (\ref{eq:A}) with
some integral conditions. A model situation is the following form of the Phragm\'en--Lindel\"of theorem.

\medskip
\noindent
{\bf Theorem E} (Ahlfors \cite{Ah}) {\em If a subharmonic function $u$ in $\C_+$ with nonpositive boundary values on $\R$ satisfies
$$\lim_{r\to\infty} r^{-1}\int_0^\pi u^+(re^{i\theta})\sin\theta\,d\theta=0,$$
then $u\le 0$ in $\C_+$.}

We will show that all the above
theorems are particular cases of results on the class
$\mathcal A$ defined below and that the $\log\log$-conditions appear as conditions for continuity of certain logarithmic potentials.

\begin{definition}\label{def:A}
Let $\nu$ be a probability measure on a segment $[a,b]$; we will
identify it occasionally with its distribution function
$\nu(t)=\nu([a,t])$. Suppose $\nu(t)$ is strictly increasing and continuous
on $[a,b]$, and denote by $\mu$ its inverse function extended to the
whole real axis as $\mu(t)=a$ for $t<0$ and $\mu(t)=b$ for $t>1$.
We will say that such a measure $\nu$ belongs to the class $\cA[a,b]$
if
\begin{equation}
\label{eq:cA}
\lim_{\delta\to 0}\sup_x\int_0^\delta\frac{\mu(x+t)-\mu(x-t)}{t}\,dt=0.
\end{equation}
\end{definition}

Note that this class is completely different from MacLane's class
$\mathcal A$ \cite{Mac} that consists of holomorphic functions in
the unit disk with asymptotic values at a dense subset of the
circle. MacLane's class is however described by the condition
$|f(re^{i\theta})|\le M(r)$, $M\in\LL[0,1]$.

Our results extending Theorems A--C and E are as follows.

\begin{theo}
\label{theo:1} Let a subharmonic function $u$ in the complex plane
satisfy
\begin{equation}
\int_0^{2\pi} u^+(te^{i\theta})\,d\nu(\theta) \le V(t)\quad \forall t\ge
t_0, \label{eq:1}
\end{equation}
with $\nu\in\cA[0,2\pi]$ and a nondecreasing function $V$ on $\R_+$. Then
there exist constants $c>0$ and $A\ge 1$, independent of $u$, such
that
\begin{equation}
\label{eq:distort} u(te^{i\theta})\le cV(At)\quad  \forall t\ge t_0.
\end{equation}
\end{theo}

\begin{theo}
\label{theo:2}
If a subharmonic function $u$
in the upper half-plane $\C_+$ satisfies the conditions
$$
\limsup_{z\to x_0} u(z)\le 0\quad\forall x_0\in\R
$$
and
$$
\lim_{t\to\infty}t^{-1}\int_0^{\pi} u^+(te^{i\theta})\,d\nu(\theta)=0
$$
with $\nu\in \cA[0,\pi]$, then $u(z)\le 0 $ $\forall z\in\C_+$.
\end{theo}

\begin{theo}
\label{theo:3} Let a subharmonic function $u$ in $Q=\{x+iy:\:
|x|<1,\ |y|<1\}$ satisfy
\begin{equation}
\int_{-1}^{1} u^+(x+iy)\,d\nu(y)\le 1\quad\forall x\in(-1,1)
\label{eq:3}
\end{equation}
with $\nu\in\cA[-1,1]$. Then for each compact set $K\subset Q$ there
is a constant $C_K$, independent of the function $u$, such that
$u(z)\le C_K$ on $ K$.
\end{theo}

Relation of these results to the $\log\log$-theorems becomes clear by means of
the following statement.

\begin{definition} Denote by $\LM[a,b]$ the class of all nonnegative
integrable functions $g$ on the segment $[a,b]$, such that
\begin{equation}
\label{eq:LM}
\int_a^b\log^-g(s)\,ds<\infty.
\end{equation}
\end{definition}

\begin{proposition}
\label{prop:1} If the density $\nu'$ of an absolutely continuous
increasing function $\nu$
 belongs to the class
$\LM[a,b]$, then $\nu\in\cA[a,b]$. Consequently, if a holomorphic function $f$ has a majorant
$M\in\LL$, then $\log|f|$ has the corresponding
integral bound with the weight $\nu\in\cA$ with the density $\nu'(t)=\min\{1,1/M(t)\}$.
\end{proposition}

We recall that positive measures $\nu$ on the unit circle with
$\nu'\in\LM[0,2\pi]$ are called {\it Szeg\"o measures}.
Proposition~\ref{prop:1} states, in particular,  that absolutely
continuous Szeg\"o measures belong to the class $\cA[0,2\pi]$.

An integral version of Theorem D has the following form.

\begin{theo}
\label{theo:4}
Let a function $u$, subharmonic in $\C$ and harmonic in $\C\setminus\R$,
satisfy the inequality
\begin{equation}
\label{eq:M1}
\int_{-\pi}^\pi u^-(re^{i\theta})\Phi(|\sin\theta|)\,d\theta\le V(r)\quad
\forall r\ge r_0,
\end{equation}
where $\Phi\in\LM[0,1]$ is nondecreasing
 and the function $V$ is such that $r^{-1-\delta}
V(r)$ is increasing in $r$ for some $\delta>0$.
Then there are constants $c>0$ and $A\ge 1$, independent of $u$,
such that
$$
u(re^{i\theta})\le cV(Ar)\quad\forall r\ge r_1=r_1(u).
$$
\end{theo}

Our proofs of Theorems \ref{theo:1}--\ref{theo:4} rest on a
presentation of measures of the class $\cA[0,2\pi]$ as radial
projections of harmonic measures of star-shaped domains. Let
$\Omega$ be a bounded Jordan domain containing the origin. Given a
set $E\subset\partial\Omega$, $\omega(z,E,\Omega)$ will denote the
harmonic measure of $E$ at $z\in\Omega$, i.e., the solution of the
Dirichlet problem in $\Omega$ with the boundary data $1$ on $E$ and $0$
on $\partial\Omega\setminus E$. The measure $\omz$ generates a
measure on the unit circle $\T$ by means of the radial projection
$\zeta\mapsto \zeta/|\zeta|$. It is convenient for us to consider it
as a measure on the segment $[0,2\pi]$, so we put
\begin{equation}
\label{eq:hat}
\widehat\omega_\Omega(F)=\omega(0,\{\zeta\in\partial\Omega:\arg\zeta\in
F\}, \Omega)
\end{equation}
for each Borel set $F\subset[0,2\pi]$.

The inverse problem is as follows. {\sl Given a probability measure on
the unit circle $\T$, is it the radial projection of the harmonic
measure of any domain $\Omega$?}

For our purposes we specify $\Omega$ to be {\it strictly star-shaped},
i.e., of the form
\begin{equation}
\label{eq:domain}
\Omega=\{re^{i\theta}:\: r<r_\Omega(\theta),\ 0\le\theta\le 2\pi\}
\end{equation}
with $r_\Omega$ a positive continuous function on $[0,2\pi]$,
$r_\Omega(0)=r_\Omega(2\pi)$.

\begin{theo}
\label{theo:5}
A continuous probability measure $\nu$ on $[0,2\pi]$ is the radial projection
of the harmonic measure of a strictly star-shaped domain if and only if
$\nu\in\cA[0,2\pi]$.
\end{theo}

\begin{corollary}
\label{cor:Szego} Every absolutely continuous measure from the
Szeg\"o class on the unit circle is the radial projection of the
harmonic measure of some strictly star-shaped domain.
\end{corollary}

Theorem \ref{theo:5} is proved by a method originated by B.Ya. Levin
in theory of majorants in classes of subharmonic functions
\cite{Levin}.

Theorems \ref{theo:1}--\ref{theo:3} and \ref{theo:5} (some of them
in a slightly weaker form) were announced in \cite{R1} and proved in
\cite{R2} and \cite{R3}. The main objective of the present paper,
Theorem~\ref{theo:4}, is new. Since its proof rests heavily on
Theorem~\ref{theo:5}, we present a proof of the latter as well,
having in mind that the papers \cite{R2} and \cite{R3} are not
easily accessible. Moreover, we include the proofs of Theorems
\ref{theo:1}--\ref{theo:3}, too, motivated by the same accessability
reason as well as by the idea of showing the whole picture.


\section{Radial projections of harmonic measures (Proofs of Theorem \ref{theo:5}
and Proposition \ref{prop:1})}

Measures from the class $\cA$ have a simple characterization as follows.

\begin{proposition} Let $\mu$ and $\nu$ be as in Definition~\ref{def:A}. Then the function $$N(x)=\int_0^1\log|x-t|\,d\mu(t)$$ is continuous on $[0,1]$ if and only if $\nu\in\cA[a,b]$.
\end{proposition}

{\it Proof.} The function $N(x)$ is continuous on $[0,1]$ if and only if for any $\epsilon>0$ one can choose $\delta\in(0,1)$ such that
$$I_x(\delta)=\int_{|t-x|<\delta}\log|x-t|\,d\mu(t)>-\epsilon$$
for all $x\in[0,1]$. Integrating $I_x$ by parts, we get
$$
|I_x(\delta)|=\int_0^\delta\frac{r_x(t)}{t}\,dt+r_x(\delta)|\log\delta|,
$$
where $r_x(t)=\mu(x+t)-\mu(x-t)$. Therefore, continuity of $N(x)$ implies (\ref{eq:cA}).
On the other hand, since $r_x(t)$ increases in $t$, we have
 $$r_x(\delta)|\log\delta|=2r_x(\delta) \int_\delta^{\sqrt\delta}\frac{dt}{t}\le 2\int_\delta^{\sqrt\delta}\frac{r_x(t)}{t}\,dt,
 $$
which gives the reverse implication. {\hfill$\square$\rm}

\medskip
In the proof of Theorem \ref{theo:5}, we will use this property in the following form.

\begin{proposition}\label{prop:cont} Let $\mu$ and $\nu$ be as in Definition~\ref{def:A} for the class $\cA[0,2\pi]$. Then the function $$h(z)=\int_0^{2\pi}\log|e^{i\theta}-z|\,d\mu(\theta/2\pi)$$ is continuous on $\T$ if and only if $\nu\in\cA[0,2\pi]$.
\end{proposition}

\medskip

{\it Proof of Theorem \ref{theo:5}}. 1) First we prove the
sufficiency: {\sl every $\nu\in\cA[0,2\pi]$ has the form
$\nu=\widehat\omega_\Omega$ (\ref{eq:hat}) for some strictly
star-shaped domain $\Omega$. In particular, for any compact set
$K\in \Omega$ there is a constant $C(K)$ such that
\begin{equation}
\label{eq:bound}
\omega(z, E,\Omega)\le C(K)\,\nu(\arg E)\quad\forall z\in E
\end{equation}
for every Borel set $E\subset\partial\Omega$, where $\arg
E=\{\arg\zeta:\: \zeta\in E\}$.}

Let
$$
u(z)={\frac1\pi}\int_0^{2\pi}\log|e^{i\theta}-z|\,d\mu({\theta /
2\pi})
$$
with $\mu$ the inverse function to $\nu\in\cA[0,2\pi]$. The function $u$ is
subharmonic in $\C$ and harmonic outside the unit circle $\T$. By Proposition~\ref{prop:cont},
it is continuous  on $\T$ and thus, by Evans' theorem, in the whole plane. Let $v$ be a
harmonic conjugate to $u$ in the unit disk $\D$, which is determined
uniquely up to a constant. Since $u\in C({\overline \D})$, radial
limits $v^*(e^{i\psi})$ of $v$ exist a.e.\ on $\T$. Let us fix such a
point $e^{i\psi_0}$ and choose the constant in the definition of $v$
in such a way that $v^*(e^{i\psi_0})=\psi_0$.

Consider then the function $w(z)=z\exp\{-u(z)-iv(z)\}$, $z\in \D$.
By the Cauchy-Riemann condition, $\partial v/\partial\phi=r\partial
u/\partial r$, which implies
\begin{eqnarray*}
\arg w(re^{i\psi}) &=& \psi-v(re^{i\psi_0})
  -\int_{\psi_0}^\psi \frac{\partial v(re^{i\phi})}{\partial\phi}\,d\phi
= \psi_0-v(re^{i\psi_0})\\
 {}& +&\frac1{2\pi}\int_{\psi_0}^\psi \int_0^{2\pi}
\left[1-\frac{2r^2-2r\cos(\theta-\phi)}{|r-e^{i(\theta-\phi)}|^2}\right]
\,d\mu({\theta/ 2\pi})\,d\phi\\
{}&=& \psi_0-v(re^{i\psi_0})
  +\frac1{2\pi}\int_{\psi_0}^\psi \int_0^{2\pi}
\frac{1-r^2}{|r-e^{i(\theta-\phi)}|^2} \,d\mu({\theta/
2\pi})\,d\phi.
\end{eqnarray*}
By changing the integration order and passing to the limit as $r\to 1$,
we derive that for each $\psi\in [0,2\pi]$ there exists the limit
$$
\lim_{r\to 1}\arg w(re^{i\psi}) =\mu({\psi/ 2\pi})
-\mu({\psi_0/ 2\pi}).
$$
Therefore the function $\arg w$ is continuous up to the boundary of
the disk; in particular, we can take $\psi_0=0$. Since $|w|$ is
continuous in $\overline \D$ as well, so is $w$.

By the boundary correspondence principle, $w$ gives a conformal map
of $\D$ onto the domain
\begin{equation}
\label{eq:obl1}
\Omega=\{re^{i\theta}:\: r<\exp\{-u(\exp\{2\pi i\nu(\theta)\})\},\
0\le\theta\le 2\pi\}.
\end{equation}
It is easy to see that the domain $\Omega$ is what we sought. Let
$f$ be the conformal map of $\Omega$ to $\D$, inverse to $w$. For
$z\in\Omega$ and $E\subset\partial\Omega$, we have
\begin{eqnarray*}
\omega(z,E,\Omega) &=& \omega(f(z),f(E),U)=\frac1{2\pi}
\int_{\arg{f(E)} }\frac{1-|f(z)|^2}{|f(z)-e^{it}|^2}\,dt\\
 {}&=& (1-|f(z)|^2)\int_{\arg E}\frac{d\nu(s)}{|f(z)-e^{2\pi i\nu(s)}|^2},
\end{eqnarray*}
which proves the claim.

\medskip
2) Now we prove the necessity: {\sl if $\omega$ is of the form
(\ref{eq:domain}), then $\widehat\omega_\Omega \in\cA[0,2\pi]$}.

We use an idea from the proof of \cite[Theorem 2.4]{Levin}.
Let $w$ be a conformal map of $\D$ to $\Omega$, $w(0)=0$.
Since $\Omega$ is a Jordan domain, $w$ extends to a continuous map
from $\overline \D$ to $\overline\Omega$, and we can specify it to
have $\arg w(1)=0$. Define
$$
f(z)=u(z)+iv(z)=\log\frac{w(z)}z\ {\rm for}\ |z|\le 1, \quad
f(z)=f(|z|^{-2}z)\ {\rm for}\ |z|>1.
$$
It is analytic in $\D$ and continuous in $\C$.
Define then the function
\begin{equation}\label{eq:lam}\lambda(z)=u(z)+\frac1{\pi}\int_0^{2\pi}\log|e^{i\psi}-z|\,dv(e^{i\psi}),
\end{equation}
$\delta$-subharmonic in $\C$ and harmonic in $\C\setminus\T$. Let us show that it as actually harmonic (and, hence, continuous) everywhere. To this end, take any function $\alpha\in C(\T)$ and a number $r<1$, and apply Green's formula for $u(z)$ and $A(z)=|z|\alpha(z/|z|)$ in the domain $D_r=\{r<|z|<r^{-1}\}$:
\begin{equation}\label{eq:12}
\int_{D_r}(A\Delta u-u\Delta A) = \left[\frac{\rho}{2\pi}\int_0^{2\pi}\left(\rho\alpha(e^{i\psi}) \frac{\partial u(\rho e^{i\psi})}{\partial \rho}-
u(\rho e^{i\psi})\alpha(e^{i\psi})\right)\,d\psi
\right]_{\rho=r}^{\rho=R}.
\end{equation}
Using the definition of the function $f$ outside $\D$ and the Cauchy-Riemann equations
$\partial v/\partial\phi=\rho\partial u/\partial \rho$ if $\rho<1$ and $\partial v/\partial\phi=-\rho\partial u/\partial \rho$ if $\rho>1$ (which follows from the definition of $f$), we can write the right hand side of (\ref{eq:12}) as
$$
-\frac{r+r^{-1}}{2\pi}\int_0^{2\pi}\alpha(e^{i\psi})\,d_\psi v(r e^{i\psi})+ \frac{r-r^{-1}}{2\pi}\int_0^{2\pi}u(r e^{i\psi})\alpha(e^{i\psi})\,d\psi.
$$
When $r\to 1$, (\ref{eq:12}) takes the form
$$\int_\T \alpha\Delta u=-\frac{1}{\pi}\int_0^{2\pi}\alpha(e^{i\psi})\,d v(e^{i\psi}),$$
which implies the harmonicity of the function $\lambda(z)$ (\ref{eq:lam}) in the whole plane.

Now we recall that $v(e^{i\psi})=\arg w(e^{i\psi})-\psi$. Since the
harmonic measure of the $w$-image of the arc
$\{e^{i\theta}:0<\theta<\psi\}$ equals $\psi/2\pi$, we have
$$
\widehat\omega_\Omega(\arg w(e^{i\psi}))=\psi/2\pi
$$
and thus $\arg w(e^{i\psi})=\mu(\psi/2\pi)$ with $\mu$ the inverse
function to $\widehat\omega_\Omega(\psi)$. Therefore,
$v(e^{i\psi})=\mu(\psi/2\pi)-\psi$.

Consider, finally, the function
$$
\gamma(z)=\frac1{\pi}\int_0^{2\pi}\log|e^{i\psi}-z|\,d\mu(\psi/2\pi)=
\lambda(z)-u(z)+
\frac1{\pi}\int_0^{2\pi}\log|e^{i\psi}-z|\,d\psi.
$$
Since it is continuous on $\T$,
 Proposition~\ref{prop:cont} implies $\widehat\omega_\Omega\in \cA[0,2\pi]$, and the theorem
is proved. {\hfill$\square$\rm}

\medskip
Note that all the dilations $t\Omega$ of $\Omega$ ($t>0$) represent
the same measure from $\cA[0,2\pi]$, and $\Omega$ with a given
projection $\widehat\omega_\Omega$ is unique up to the dilations.

\medskip
Now we prove Proposition \ref{prop:1} that presents a wide subclass
of $\cA$ with a more explicit description.

\medskip

{\it Proof of Proposition \ref{prop:1}}. Let $\nu:[0,1]\to [0,1]$
be an absolutely continuous, strictly increasing function, $\nu'\in\LM[0,1]$.
Since ${\rm mes}\,\{t:\nu'(t)=0\}=0$, its inverse function $\mu$ is absolutely
continuous (\cite{N}, p. 297), so
$$
\mu(t)=\int_0^t g(s)\,ds
$$
with $g$ a nonnegative function on $[0,1]$. We have
$$
\infty>\int_0^1\log^-\nu'(t)\,dt=\int_0^1\log^-\frac1{\mu'(t)}\,d\mu(t)
=\int_0^1g(t)\log^+g(t)\,dt,
$$
so $g$ belongs to the Zygmund class ${\bf L\log L}$.

Let $\Delta(t) $ denote the modulus of continuity of the function $\mu$.
Note that it can be expressed in the form
$$
\Delta(t)=\int_0^t h(s)\,ds
$$
where $h$ is the nonincreasing equimeasurable rearrangement of $g$. Then
\begin{eqnarray*}
\int_0^1\frac{\Delta(t)}{t}\,dt &=&\int_0^1 t^{-1}\int_0^1
h(s)\,ds\,dt
=\int_0^1h(s)\log s^{-1}\,ds\\
{}&=& \int_{E_1\cup E_2}h(s)\log s^{-1}\,ds,
\end{eqnarray*}
where $E_1=\{s\in (0,1):h(s)>s^{-1/2}\}$, $E_2=(0,1)\setminus E_1$.
Since $h\in {\bf L\log L}\,[0,1]$,
$$
 \int_{E_1}h(s)\log s^{-1}\,ds\le 2 \int_{E_1}h(s)\log h(s)\,ds
<\infty.
$$
Besides,
$$ \int_{E_2}h(s)\log s^{-1}\,ds\le  \int_{E_2}s^{-1/2}\log s^{-1}\,ds<\infty.
$$
Therefore,
$$
\int_0^1\frac{\Delta(t)}{t}\,dt<\infty
$$
and thus
$$
\lim_{\delta\to 0}\int_0^\delta \frac{\Delta(t)}{t}\,dt=0,
$$
which gives (\ref{eq:cA}). {\hfill$\square$\rm}

\medskip
Corollary \ref{cor:Szego} follows directly from the definition of the Szeg\"o
class, Theorem~\ref{theo:5} and Proposition~\ref{prop:1}.


\section{Proofs of Theorems \ref{theo:1} and \ref{theo:2}}

Here we show how the integral variants of Carleman's and Wolf's theorems
can be derived from Theorem \ref{theo:5}.

We will need an elementary

\begin{lemma}
\label{lem:1} Let $r(\theta)\in C[0,2\pi]$, $1<r_1\le r(\theta)\le
r_2$, let $\nu$ be a positive measure on $[0,2\pi]$ and $V(t)$ be a
nonnegative function on $[0,\infty]$. If a nonnegative function
$v(te^{i\theta})$ satisfies
$$
\int_0^{2\pi}v(te^{i\theta})\,d\nu(\theta)\le V(t)\quad\forall t\ge  t_0,
$$
then for any $R_2>R_1\ge t_0$,
$$
\int_{R_1}^{R_2}\int_0^{2\pi}v(t\,r(\theta)e^{i\theta})\,d\nu(\theta)\,dt
\le r_1^{-1}\int_{r_1R_1}^{r_2R_2} V(t)\,dt.
$$
\end{lemma}

{\it Proof of Lemma \ref{lem:1}} is straightforward:
$$
\int_{R_1}^{R_2}\int_0^{2\pi}v(t\,r(\theta)e^{i\theta})\,d\nu(\theta)\,dt
=\int_0^{2\pi}\int_{R_1r(\theta)}^{R_2r(\theta)}v(te^{i\theta})
\,dt\,\frac{d\nu(\theta)}{r(\theta)}
$$
$$
\le r_1^{-1}\int_0^{2\pi}\int_{R_1r_1}^{R_2r_2}v(te^{i\theta})
\,dt\,{d\nu(\theta)}
\le r_1^{-1}\int_{R_1r_1}^{R_2r_2}V(t)
\,dt.
$$
{\hfill$\square$\rm}

\medskip
{\it Proof of Theorem \ref{theo:1}}. By Theorem \ref{theo:5}, there
exists a domain $\Omega$ of the form (\ref{eq:domain}) that contains
$\overline \D$ such that
\begin{equation}
\label{eq:t11} \omega(z,E,\Omega)\le c_1\nu(\arg E),\quad
\forall z\in{\overline \D}, \  E\subset\partial \Omega,
\end{equation}
with a constant $c_1>0$, see (\ref{eq:bound}). Let $r_1=\min
r(\theta)$. $r_2=\max r(\theta)$. By the Poisson--Jensen formula
applied to the function $v_t(z)=u^+(tz)$ ($t>0$) in the domain
$s\Omega$ ($s>1$) we have, due to (\ref{eq:t11}),
\begin{eqnarray*}
v_t(z) &\le &\int_{\partial s\Omega}
v_t(\zeta)\,\omega(z,d\zeta,s\Omega)
=\int_{\partial \Omega} v_t(s\zeta)\,\omega(s^{-1}z,d\zeta,\Omega)\\
{}&\le & c_1\int_0^{2\pi}v_t(s\,r(\theta)e^{i\theta})\,d\nu(\theta),
\quad z\in{\overline \D}.
\end{eqnarray*}
The integration of this relation over $s\in [1,R]$ ($R>1$) gives,
by Lemma \ref{lem:1},
$$
(R-1)v_t(z)\le  c_1\int_1^R\int_0^{2\pi}v_t(s\,r(\theta)
e^{i\theta})\,d\nu(\theta)\,ds\le
c_2t^{-1}r_1^{-1}\int_{tr_1}^{tr_2R}V(s)\,ds
$$
for each $t\ge t_0$.
So,
$$
u(te^{i\theta})\le c(R)V(t\,r_2R),\quad t\ge t_0,
$$
which proves the theorem.{\hfill$\square$\rm}

\medskip
{\bf Remarks}. 1. It is easy to see that the constant $A$ in (\ref{eq:distort})
can be chosen arbitrarily close to $r_2/ r_1\ge 1$.

2. Note that we have used inequality (\ref{eq:1}) in the integrated
form only, so the following statement is actually true: {\em If a
subharmonic function $u$ on $\C$ satisfies
\begin{equation}
\int_ {t_0}^t\int_0^{2\pi} u^+(se^{i\theta})\,d\nu(\theta)\,ds \le
W(t)\quad\forall t\ge t_0 \label{eq:1'}
\end{equation}
with $\nu\in\cA[0,2\pi]$ and a nondecreasing function $W$, then
there are constants $c>0$ and $A\ge 1$, independent of $u$, such
that $ u(te^{i\theta})\le ct^{-1}W(At)$ for all $t\ge t_0$.}

\medskip
Now we prove Theorem \ref{theo:2} as a consequence of Theorem \ref{theo:1}.

\medskip
{\it Proof of Theorem \ref{theo:2}}. The function $v$ equal to $u^+$
in $\C_+$ and $0$ in $\C\setminus\C_+$ is a subharmonic function in
$\C$ satisfying the condition
$$
\int_0^{2\pi} v^+(te^{i\theta})\,d\nu(\theta)\le V_1(t)
$$
with $\nu\in \cA[0,2\pi]$ and $V_1(t)=o(t)$, $t\to\infty$.
Therefore, it satisfies the conditions of Theorem \ref{theo:1} with
the majorant $V(t)=\sup\{V_1(s):s\le t\}$. So, $ \sup_\theta
u^+(te^{i\theta})=o(t)$ as $t\to\infty$, and the conclusion holds by
the standard Phragm\'en--Lindel\"of theorem. {\hfill$\square$\rm}


\section{Proof of Theorem \ref{theo:3}}

The integral version of the Levinson--Sj{\"o}berg theorem will be proved
along the same lines as Theorem \ref{theo:1}, however the local situation
needs a more refined adaptation.

We start with two elementary statements close to Lemma \ref{lem:1}.

\begin{lemma}
\label{lem:2} Let a nonnegative integrable function $v$ in the
square $Q=\{|x|,|y|<1\}$ satisfy (\ref{eq:3}) with a continuous strictly
increasing function $\nu$. Then for any $d\in(0,1)$ there exists a
constant $M_1(d)$, independent of $u$, such that for each
$y_0\in(-1,1)$ one can find a point $y_1\in(-1,1)\cap (y_0-d,y_0+d)$
with
$$
\int_{-1}^1v(x+iy_1)\,dx<M_1(d).
$$
\end{lemma}

{\it Proof}. Assume $y_0\ge 0$, then
$$
\int_{y_0-d}^{y_0} \int_{-1}^1v(x+iy)\,dx \,d\nu(y)=
\int_{-1}^1\int_{y_0-d}^{y_0} v(x+iy)\,d\nu(y)\,dx\le 2.
$$
Therefore for some $y_1\in(y_0-d,y_0)$,
$$
\int_{-1}^1v(x+iy_1)\,dx\le 2[\nu(y_0)-\nu(y_0-d)]^{-1}\le
2[\Delta_*(\nu,d)]^{-1}
$$
with $\Delta_*(\nu,d)=\inf\{\nu(t)-\nu(t-d):\,t\in(0,1)\}>0$.
{\hfill$\square$\rm}

\begin{lemma}
\label{lem:3}
 Let a function $v$ satisfy the conditions of Lemma \ref{lem:2}, a
 function $r$ be continuous on a segment $[a,b]\subset[-1,1]$,
$0<r_1=\min r(y)\le\max r(y)=r_2<1$,
 and $\delta\in (0,1-r_2)$.
Then there exists $t\in(0,\delta)$ such that
$$
\int_a^b v(t+r(y)+iy)\,d\nu(y)<M_2(\delta)
$$
 with  $M_2(\delta)$
independent of $v$.
\end{lemma}

{\it Proof}. We have
\begin{eqnarray*}
\int_0^\delta\int_a^b v(t+r(y)+iy)\,d\nu(y) &=&
\int_a^b\int_{r(y)}^{\delta+r(y)} v(s+iy)\,ds\,d\nu(y)\\
\le\int_{r_1}^{\delta+r_2}\int_a^b v(s+iy)\,d\nu(y)\,ds &\le &
\delta+r_2-r_1.
\end{eqnarray*}
Thus one can find some $t\in(0,\delta)$ such that
$$
\int_a^b v(t+r(y)+iy)\,d\nu(y)<\delta^{-1}(\delta+r_2-r_1).
$$
{\hfill$\square$\rm}

\medskip
{\it Proof of Theorem \ref{theo:3}}. Consider the measure $\nu_1$ on
$[-i,i]$ defined as $$\nu_1(E)=\nu(-iE),\quad E\subset [-i,i].$$ The
conformal map $f(z)=\exp\{z\pi/2\}$ of the strip $\{|\Im
z|<1\}$ to the right half-plane $\C_r$ pushes the measure $\nu_1$ forward to the
measure $f^*\nu$ on the semicircle $\{e^{i\theta}: -\pi/2\le
\theta\le\pi/2\}$, producing a measure of the
class $\cA[-\pi/2,\pi/2]$; we extend it to some measure
$\nu_2\in\cA[-\pi.\pi]$. By Theorem~\ref{theo:5}, there is a strictly
star-shaped domain $\Omega\supset{\overline \D}$ such that the
radial projection of its harmonic measure at $0$ is the
normalization $\nu_2/\nu_2([-\pi,\pi])$ of $\nu_2$.

Let $\Omega_1=\Omega\cap\C_r$, then for every Borel set $E\subset\Gamma=\partial
\Omega_1\cap\C_r$ and any compact set $K\subset\Omega_1$,
$$
\omega(w,E,\Omega_1)\le C_1(K)\,\nu_2(\arg E)\quad\forall w\in K.
$$

The pre-image $\Omega_2=f^{[-1]}(\Omega_1)$ of $\Omega_1$ has the form
$$
\Omega_2=\{z=x+iy:\: x<\varphi(y),\ y\in (0,1)\}
$$
with some function $\varphi\in C[-1,1]$. Let
$$
\Gamma_2=\{x+iy:\: x=\varphi(y),\ y\in (0,1)\},
$$
then for every Borel $E\subset\Gamma_2$ and any compact subset $K$ of $\Omega_2$,
\begin{equation}
\label{eq:LS1} \omega(z,E,\Omega_2)\le C_2(K)\,\nu(\Im E)\quad \forall z\in  K.
\end{equation}
For the domain
$$
\Omega_3=\{z=x+iy:\: x>-\varphi(y),\ y\in (0,1)\}
$$
we have, similarly, the relation
\begin{equation}
\label{eq:LS2} \omega(z,E,\Omega_3)\le C_3(K)\,\nu(\Im E)\quad \forall z\in
K
\end{equation}
for each $E\subset\Gamma_3=\{x+iy:\: x=-\varphi(y),\ y\in (0,1)\}$ and compact set $K\subset\Omega_3$.

Let now $K$ be an arbitrary compact subset of the square $Q$. We
would be almost done if we were able to find some reals $h_2(K)$ and
$h_3(K)$ such that
$$K\subset\{\Omega_2+h_2(K)\}\cap\{\Omega_3+h_3(K)\}\subset\overline{\{\Omega_2+h_2(K)\}\cap\{\Omega_3+h_3(K)\}} \subset Q.$$
However this is not the case for any $K$ unless $\varphi\equiv
const$. That is why we need partition.

Given $K$ compactly supported in $Q$, choose a positive $\lambda<(4\, {\rm
dist}\,(K,\partial Q))^{-1}$ and then $\tau\in(0,\lambda)$ such that
the modulus of continuity of $\varphi$ at $4\tau$ is less than
$\lambda$. Take a finite covering of $K$ by disks
$B_j=\{z:|z-z_j|<\tau\}$, $z_j\in K$, $1\le j\le n$. To prove the
theorem, it suffices to estimate the function $u$ on each $B_j$.

Let $Q_j=\{z\in Q: |\Im (z-z_j)|<2\tau\}$, then $B_j\subset
Q_j$ and ${\rm dist}\,(B_j,\partial Q_j)=\tau$. Take also
$$
\Omega_2^{(j)}=\Omega_2\cap Q_j,\quad \Gamma_2^{(j)}=
\Gamma_2\cap{\overline\Omega_2^{(j)}}=\{x+iy:\: x=\varphi(y),\
a_j\le y\le b_j\}.
$$
Now we can find reals $h_2^{(j)}$ and $h_3^{(j)}$ such that
$$
\Gamma_2^{(j)} +h_2^{(j)}
=\{x+iy: x=r_2^{(j)}(y)\}
\subset Q_j\cap\{x+iy: 1-4\lambda <x<1<2\lambda\}
$$
and
$$
\Gamma_3^{(j)} +h_3^{(j)}
=\{x+iy: x=r_3^{(j)}(y)\}
\subset Q_j\cap\{x+iy: -1+2\lambda<x<-1+4\lambda\}.
$$
Furthermore, by Lemma \ref{lem:3}, there exist $t_2^{(j)}\in(0,\lambda)$
and $t_3^{(j)}\in(-\lambda,0)$
such that
\begin{equation}
\label{eq:LS18}
\int_{a_j}^{b_j} u^+(t_k^{(j)}+r_k^{(j)}(y)+iy)\,d\nu(y)<M_2(\lambda),
\quad k=2,3.
\end{equation}
Finally we can find, due to Lemma \ref{lem:2},
$y_1^{(j)}\in (a_j,a_j+\tau)$ and $y_2^{(j)}\in (b_j-\tau,b_j)$
such that
\begin{equation}
\label{eq:LS20}
\int_{-1}^1u^+(x+iy_m)\,dx<M_1(\tau),\quad m=1,2.
\end{equation}
Denote
$$
\Omega^{(j)}=\{x+iy:r_3^{(j)}(y)+t_3^{(j)}<x<
r_2^{(j)}(y)+t_2^{(j)},\ y_1^{(j)}\le y\le y_2^{(j)}\}.
$$
Since $\overline{B_j}\subset\Omega^{(j)}$, relations (\ref{eq:LS1}) and (\ref{eq:LS2})
imply
\begin{equation}
\label{eq:LS21}
\omega(z,E,\Omega^{(j)})\le C(B_j)\nu(\Im E)\quad\forall z\in B_j
\end{equation}
for all $E$ in the vertical parts of $\partial\Omega^{(j)}$. For $E$
in the horizontal parts of $\partial\Omega^{(j)}$, we have,
evidently,
\begin{equation}
\label{eq:LS22}
\omega(z,E,\Omega^{(j)})\le C(B_j)\,{\rm mes}\,E\quad\forall z\in B_j.
\end{equation}

Now we can estimate $u(z)$ for $z\in B_j$. By
(\ref{eq:LS18})--(\ref{eq:LS22}),
\begin{eqnarray*}
u(z) &\le & \int_{\partial\Omega^{(j)}}u^+(\zeta)\omega(z,d\zeta,
\Omega^{(j)})\\
{}&\le & C(B_j)\sum_{k=2}^3\int_{a_j}^{b_j}u^+(t_k^{(j)}+r^{(j)}(y)
+iy)\,d\nu(y)\\
{}& + & C(B_j)\sum_{m=1}^2\int_{-1}^1 u^+(x+iy_m)\,dx\\
{}&\le & 2C(B_j)(M_1(\tau)+M_2(\lambda)),
\end{eqnarray*}
which completes the proof. {\hfill$\square$\rm}


\section{Proof of Theorem \ref{theo:4}}
By Theorem \ref{theo:1} and Proposition \ref{prop:1}, it suffices to prove

\begin{proposition}
\label{prop:2}
If a function $u$ satisfies the conditions of Theorem \ref{theo:4},
then there exists a function $f\in\LM[-\pi,\pi]$ and
a constant $c_1>0$, the both independent of $u$,  such that
\begin{equation}
\label{eq:t41}
\int_{-\pi}^{\pi}u^+(re^{i\theta})f(\theta)\,d\theta\le
c_1V(r)\quad\forall  r>r_0.
\end{equation}
\end{proposition}

{\it Proof}. What we will do is a refinement of the arguments from the
proof of the original Matsaev's theorem (see
\cite{M}, \cite{Levin1}). Let
$$
D_{r,R,a}=\{z\in\C:\: r<|z|<R,\ |\arg z-{\pi/ 2}|<\pi({1/ 2}-a)\},
\quad 0<a<1/4,
$$
$b=(1-2a)^{-1}$, $S(\theta,a)=\sin b(\theta-a\pi)$. Carleman's formula
for the function $u$ harmonic in $D_{r,R,a}$ has the form
$$
2bR^{-b}\int_{\pi a}^{\pi-\pi a} u(Re^{i\theta})S(\theta,a)\,d\theta
-b(r^{-b}+r^bR^{-2b})\int_{\pi a}^{\pi-\pi a} u(re^{i\theta})S(\theta,a)\,d\theta
$$
$$
-(r^{-b+1}-r^{b+1}R^{-2b})\int_{-\pi a}^{\pi a} u'_r(re^{i\theta})
S(\theta,a)\,d\theta
$$
$$
+b\int_r^R\left[u(xe^{i\pi a})+u(xe^{i\pi (1-a)})\right](x^{-b-1}
-x^{b-1}R^{-2b})\,dx=0.
$$
It implies the inequality
$$
\int_{\pi a}^{\pi-\pi a} u^+(Re^{i\theta})S(\theta,a)\,d\theta\le
c(r,u)R^b+
\int_{\pi a}^{\pi-\pi a} u^-(Re^{i\theta})S(\theta,a)\,d\theta
$$
\begin{equation}
\label{eq:M3}
+R^b\int_r^R\left[u^-(xe^{i\pi a})+u^-(xe^{i\pi (1-a)})\right](x^{-b-1}
-x^{b-1}R^{-2b})\,dx.
\end{equation}
Fix some $\tau\in(0,1/4)$ such that
\begin{equation}
\label{eq:M4}
\beta:=(1-2\tau)^{-1}<1+\delta
\end{equation}
with $\delta$ as in the statement of Theorem \ref{theo:4}.
Inequality (\ref{eq:M3}) gives us the relation
\begin{eqnarray}
\label{eq:M5} I_0 &:= & \int_0^\tau\Phi(\sin\pi a) \int_{\pi a}^{\pi-\pi a} u^+(Re^{i\theta})S(\theta,a)
\,d\theta\,da\nonumber\\
{}&\le &
c(r,u)\int_0^\tau R^b\Phi(\sin\pi a)\,da+
\int_0^\tau\Phi(\sin\pi a)
 \int_{\pi a}^{\pi-\pi a} u^-(Re^{i\theta})S(\theta,a)\,d\theta\,da\nonumber\\
{}&+& \int_0^\tau \Phi(\sin\pi a)
\int_r^R\left[u^-(xe^{i\pi a})+u^-(xe^{i\pi (1-a)})\right]
R^bx^{-b-1}
\,dx\,da\nonumber\\
{}&=& I_1+I_2+I_3.
\end{eqnarray}
We can represent $I_0$ as
$$
I_0=\int_0^\pi u^+(Re^{i\theta})\Psi(\theta)\,d\theta
$$
with
\begin{equation}
\label{eq:psi}
\Psi(\theta)=\int_0^{\lambda(\theta)} S(\theta,a)\Phi(\sin\pi a)\,da
\end{equation}
and
\begin{equation}
\label{eq:lambda}
\lambda(\theta)=\min\{\theta/\pi, 1-\theta/\pi,\tau\}.
\end{equation}

Note that $S(\theta,a)\ge 0$ when $a\le\lambda(\theta)$, and
$S'_a(\theta,a)\le 0$ for all $a<1/4$. Since $\Phi(t)$ is
nondecreasing, this implies the bound
$$
\Psi(\theta) \ge \int_{\lambda(\theta)/2}^{\lambda(\theta)}
S(\theta,a) \Phi(\sin\pi a)\,da\ge f(\theta)=\lambda^2(\theta)\: \Phi\left(\sin\frac{\pi\lambda(\theta)}{2}\right)
$$
and thus,
\begin{equation}
\label{eq:M9}
I_0\ge \int_0^\pi u^+(Re^{i\theta})f(\theta)\,d\theta
\end{equation}
with $f\in\LM[0,\pi]$.

Let us now estimate the right hand side of (\ref{eq:M5}). We have
\begin{equation}
\label{eq:M11}
I_1\le c(r,u)R^\beta\int_0^\tau\Phi(\sin\pi a)\,da\le c_1(r,\tau,u)R^\beta;
\end{equation}
\begin{equation}
\label{eq:M12}
I_2=\int_0^\pi u^-(Re^{i\theta})\Psi(\theta)\,d\theta
\le \int_0^\pi u^-(Re^{i\theta})\Phi(\sin\theta)\,d\theta;
\end{equation}
\begin{eqnarray}
\label{eq:M13} I_3 &\le & \int_0^\tau\int_r^R \Phi(\sin\pi a)
\left[u^-(xe^{i\pi a})+u^-(xe^{i\pi (1-a)})\right]
\left(\frac{R}{x}\right)^{\beta}x^{-1}\,dx\,da
\nonumber\\
{}&=& R^\beta\int_r^R x^{-\beta-1}\left[\int_0^{\pi\tau}+
\int_{\pi(1-\tau)}^\pi\right] u^-(xe^{i\theta}) \Phi(\sin\theta)\,d\theta\,dx
\nonumber\\
{}&\le & R^\beta\int_r^R x^{-\beta-1}\int_0^{\pi}
u^-(xe^{i\theta}) \Phi(\sin\theta)\,d\theta\,dx.
\end{eqnarray}
We insert (\ref{eq:M9})--(\ref{eq:M13}) into (\ref{eq:M5}):
\begin{eqnarray}
\label{eq:M14}
\int_0^\pi u^+(Re^{i\theta})f(\theta)\,d\theta &\le & c_1(r,\tau,u)R^\beta
+\int_0^\pi u^-(Re^{i\theta})\Phi(\sin\theta)\,d\theta\nonumber\\
{}&+& R^\beta\int_r^R x^{-\beta-1}\int_0^{\pi}
u^-(xe^{i\theta}) \Phi(\sin\theta)\,d\theta\,dx\nonumber\\
{}&=& J_1(R)+J_2(R)+J_3(R).
\end{eqnarray}
By the choice of $\beta$ (\ref{eq:M4}), $J_1(R)=o(V(R))$ as $R\to\infty$.
Condition (\ref{eq:M1}) implies $J_2(R)\le V(R)$, $R>r_0$.
As to the term $J_3$, take any $\epsilon\in(0,1+\delta-\beta)$,
then
\begin{eqnarray*}
J_3(R) &\le & R^\beta\int_r^R x^{-\beta-1}V(x)\,dx
  = R^\beta\int_r^R x^{-\beta-\epsilon}V(x)x^{\epsilon-1}\,dx\\
{}&\le &  R^\beta  R^{-\beta-\epsilon}V(R)\int_r^R x^{\epsilon-1}\,dx
\le \epsilon^{-1}V(R).
\end{eqnarray*}
These bounds give us
$$
\int_0^\pi u^+(Re^{i\theta})f(\theta)\,d\theta\le c_2V(R)\quad \forall
R>r_1(u).
$$

Absolutely the same way, we get a similar inequality in the lower
half-plane and, as a result, relation (\ref{eq:t41}).
{\hfill$\square$\rm}

\medskip
{\it Remark.} We do not know if condition (\ref{eq:M1}) can be replaced by a more general one in terms of the class $\cA$.

\medskip
{\it Acknowledgement.} The author is grateful to Alexandre Eremenko and Misha Sodin for valuable discussions, and to the referee for suggestions that have simplified considerably the proof of Theorem~5.


\vskip1cm

Tek/Nat, University of Stavanger, 4036 Stavanger, Norway

\vskip0.1cm

{\sc E-mail}: alexander.rashkovskii@uis.no

\end{document}